\documentclass[12pt,draft]{amsart}
\usepackage{amsmath,amsthm,latexsym,amscd,amsbsy,amssymb,pb-diagram}
\chardef\bslash=`\\ % p. 424, TeXbook

\makeatletter
\def\verbatim{\interlinepenalty\@M \@verbatim
  \leftskip\@totalleftmargin\advance\leftskip2pc
  \frenchspacing\@vobeyspaces \@xverbatim}
\makeatother
\makeatletter
  \def\dgt@k{\dg@DX=-3 \dg@DY=2 \dg@SIZE=3} 
\makeatother

\makeatletter
  \def\dgt@kk{\dg@DX=3 \dg@DY=-1 \dg@SIZE=3}%
\makeatother

\theoremstyle{plain}
\newtheorem{thm}{Theorem}[section]
\newtheorem{cor}[thm]{Corollary}
\newtheorem{lem}[thm]{Lemma}
\newtheorem{pro}[thm]{Proposition}

\theoremstyle{definition}
\newtheorem{rem}[thm]{Remark}
\newtheorem{defin}[thm]{Definition}

\numberwithin{equation}{section}

\newcounter{rmnum}

\def\symbolnote#1#2{\let\thefootn=\thefootnote%
\renewcommand{\thefootnote}{\fnsymbol{footnote}}%
\footnotemark[#1]%
\footnotetext[#1]{#2}%
\let\thefootnote=\thefootn
}

\newfont{\bbb}{msbm10 scaled \magstep1}
\newfont{\bbc}{msbm8 scaled \magstep0}

\newcommand{\N}{\mbox{\bbb N}}

\begin{document}

%%%%%%% Begin Topmatter %%%%%%%%%%

\title[Bounded rank of $C^{\ast}$-algebras]
{Bounded rank of $C^{\ast}$-algebras}
\author{Alex Chigogidze}
\address{Department of Mathematics and Statistics,
University of Saskatche\-wan,
McLe\-an Hall, 106 Wiggins Road, Saskatoon, SK, S7N 5E6,
Canada}
\email{chigogid@math.usask.ca}
\thanks{The first author was partially supported by NSERC research grant.}
\author{Vesko Valov}
\address{Department of Mathematics and Computer Science, Nipissing University,
100 College Drive, P.O. Box 5002, North Bay, ON, P1B 8L7, Canada}
\email{veskov@unipissing.ca}
\thanks{The second author was partially supported by Nipissing University Research Council Grant.}
\keywords{Real rank, bounded rank, direct limit}
\subjclass{Primary: 46L05; Secondary: 54F45}

%%%%%%% End topmatter %%%%%%%%%

\begin{abstract}{We introduce a concept of the bounded rank (with respect to a positive constant) for unital $C^{\ast}$-algebras as a modification of the usual real rank and present a series of conditions insuring that bounded and real ranks coincide. These observations are then used to prove that for a given $n$ and $K > 0$ there exists a separable unital $C^{\ast}$-algebra $Z_{n}^{K}$ such that every other separable unital $C^{\ast}$-algebra of bounded rank with respect to $K$ at most $n$ is a quotient of $Z_{n}^{K}$.}
\end{abstract}

\maketitle
\markboth{A.~Chigogidze, V.~Valov}{Bounded rank of $C^{\ast}$-algebras}

\section{Introduction}\label{S:intro}
The concept of the real rank plays an important role in a variety of problems related to general classification problems of $C^{\ast}$-algebras. Despite of this widely recognized fact we still see continuing attempts of defining ``right dimension" (such as stable \cite{rieffel}, analytic \cite{murphy2}, tracial \cite{lin}, exponential \cite{phillips}, completely positive \cite{winter} ranks) for unital $C^{\ast}$-algebras. Definition of the real rank \cite{brownped91} (as well as of
its prototype -- topological stable rank  \cite{rieffel}) of unital $C^{\ast}$-algebras is based on the following standard result
from classical dimension theory (see \cite{hw}, \cite{nagata}, \cite{pears}, \cite{eng78}) characterizing the Lebesgue
dimension $\dim$ of compact spaces: 

\begin{itemize}
\item[ ]
{\em The Lebesgue dimension $\dim X$ of a compact space $X$ is the least integer $n$ such that the set $\{ f \colon C(X,\mathbb{R}^{n+1}) \colon \mathbf{0} \not\in f(X)\}$ is dense in the space\footnote{Compact-open topology is being considered.} $C(X,\mathbb{R}^{n+1})$ of all continuous maps of $X$ in to the Euclidean space $\mathbb{R}^{n+1}$.}
\end{itemize}

A map $f \in C(X,\mathbb{R}^{n+1})$ may be identified with the $(n+1)$-tuple $(\pi_{1}\circ f, \dots ,\pi_{n+1}\circ f)$, where $\pi_{k} \colon \mathbb{R}^{n+1} \to \mathbb{R}$ is a projection onto the $k$-th coordinate, $k = 1,\dots ,n+1$. The condition $\mathbf{0} = (0,\dots ,0) \not\in f(X)$ may be equivalently expressed as the condition  $\sum_{k=1}^{n+1}f_{k}^{2}(x) \neq 0$ for any $x \in X$. This, in turn, can be equivalently rephrased as the invertibility of the element $\sum_{k=1}^{n+1}f_{k}^{2}$. These two observations lead us to the definition of the real rank: 
\begin{itemize}
\item[ ]
{\em The real rank\; $\operatorname{rr}(X)$ of a unital $C^{\ast}$-algebra is the least integer $n$ such that each $(n+1)$-tuple $(x_{1},\dots ,x_{n+1})$ of self-adjoint elements of $X$ can be arbitrarily closely approximated by another $(n+1)$-tuple $(y_{1},\dots ,y_{n+1})$ of self-adjoint elements so that the element $\sum_{k=1}^{n+1}y_{k}^{2}$ is invertible.}
\end{itemize}

\noindent Of course, the analogy is quite formal and does not go far beyond the obvious observation -- $\operatorname{rr}(X) = \dim\Omega (X)$, where $\Omega (X)$ is the spectrum of the commutative unital $C^{\ast}$-algebra $X$, -- and few other straightforward extensions of certain basic facts from classical dimension theory to a non-commutative situation. But there are difficulties in finding proper algebraic interpretations in terms of the concept of the real rank of a non-commutative unital $C^{\ast}$-algebra $X$ of such an immediate geometric consequence of the condition $\dim\Omega (X) \leq n$ as the possibility not only to remove the image of any $f \colon \Omega (X) \to \mathbb{R}^{n+1}$ from $\mathbf{0}$, but even to push this image by an  $\epsilon$-move outside the open ball $O(\mathbf{0},\epsilon )$ of radius $\epsilon >0$. As it turns out in the presence of some form of functional calculus (as in the commutative or the real rank zero cases -- see Corollaries \ref{C:rrcommutative} and \ref{C:rrzero}) there exists a satisfactory analog of the above mentioned geometric fact. As for the general case, the situation remains unclear and we, as a consequence, are unable to answer the question

\begin{itemize}
\item[ ]
{\em Is it true that  $\operatorname{rr}\left( \prod\{ X_{t} \colon t \in T\right) \leq n$ for any collection of unital $C^{\ast}$-algebras $X_{t}$ such that $\operatorname{rr}(X_{t}) \leq n$ for each $t \in T$?}
\end{itemize}

The corresponding topological fact -- $\dim \beta \left(\oplus\{ \Omega (X_{t})\colon t \in T\}\right) \leq n$  -- is easy to establish\footnote{$\beta K$ stands for the Stone-\v{C}ech compactification of a space $K$.}. Perhaps the easiest way of proving the desired inequality in terms of the above given characterization of $\dim$ is first to approximate (as close as we wish) each of the restrictions $f_{t} \colon \Omega (X_{t}) \to \mathbb{R}^{n+1}$, $t \in T$, of an arbitrarily given map $f \colon  \beta \left(\oplus\{ \Omega (X_{t})\colon t \in T\}\right) \to \mathbb{R}^{n+1}$, by maps $g_{t} \colon \Omega (X_{t}) \to \mathbb{R}^{n+1}\setminus \{ \mathbf{0}\}$ and then to move each of the images $g_{t}\left(\Omega (X_{t})\right)$ outside of the open ball $O$ (centered at $\mathbf{0}$) of the appropriate radius by small moves (independent of $t \in T$) . This way we get a map $g \colon \oplus\{ \Omega (X_{t}) \colon t \in T\} \to \mathbb{R}^{n+1} \setminus O$. Since the image $g ( \oplus\{ \Omega (X_{t}) \colon t \in T\})$
has a compact closure, $g$ can be extended to a map
 $\widetilde{g} \colon \beta \left(\oplus\{ \Omega (X_{t})\colon t \in T\}\right) \to \mathbb{R}^{n+1}\setminus O$ which would be the required approximation of $f$ whose image misses $\mathbf{0}$.

Analysis of this elementary reasoning leads us to the concept of the {\em bounded} rank (to be more specific -- bounded rank with respect to an apriori given positive constant) which is defined in Section \ref{S:rjs} (we intend to extend this definition to non-unital $C^{\ast}$-algebras and study its properties in a separate note). We present certain geometric conditions guaranteeing the coincidence of the real and bounded ranks. Despite the fact that there is no satisfactory theory of joint spectra for non-commuting tuples of unital $C^{\ast}$-algebras, these conditions are formulated in terms of properties of joint spectra in order to emphasize a ``non-geometrical nature" of possible differences between the real and bounded ranks. Bounded and real ranks coincide for commutative unital $C^{\ast}$-algebras as well as for unital $C^{\ast}$-algebras of real rank zero.

One of the main advantages of the bounded rank is that the corresponding analog of the above stated question has the positive answer (Proposition \ref{P:product}). This fact is then used in the construction of separable unital universal $C^{\ast}$-algebra of the given bounded rank (Theorem \ref{T:main}). A motivation for such a result lies, once again, in the classical dimension theory. It is a well-known observation (see, for instance,
\cite[Theorem 1.3.15]{eng78}) that the Menger cube
$\mu^{n}$ contains a 
topological copy of any at most $n$-dimensional metrizable compact space.
This means that every commutative separable unital
$C^{\ast}$-algebra of real rank at most $n$ is a quotient of the $C^{\ast}$-algebra
$C\left( \mu^{n}\right)$. We extend this result to the non-commutative case. Similar result for $n = 0$ has been obtained in \cite{chi001}.

We recall the concept of a ${\mathcal C}$-invertibility introduced in \cite{chi001}. For a given class ${\mathcal C}$ of separable unital
$C^{\ast}$-algebras the
${\mathcal C}$-invertibility of a unital $\ast$-homomorphism $p \colon Y \to Z$
means that for any unital $\ast$-homomorphism
$g \colon Y \to X$, with
$X \in {\mathcal C}$, there exists a unital $\ast$-homomorphism $h \colon Z \to X$
such that $g = h\circ p$. It is easy to see that if there exists a ${\mathcal C}$-invertible unital $\ast$-homomorphism $p \colon C^{\ast}\left( {\mathbb F}_{\infty}\right)\to Z_{C}$ such that $Z_{C} \in {\mathcal C}$, 
where  $C^{\ast}\left( {\mathbb F}_{\infty}\right)$  denotes the group $C^{\ast}$-algebra of the free group on countable number of generators, then $Z_{C}$ is an universal element in the class ${\mathcal C}$. 
Indeed, since every element $X$ of ${\mathcal C}$ can be represented as the image of
$C^{\ast}\left( {\mathbb F}_{\infty}\right)$ under a
surjective $\ast$-homomorphism
$g \colon C^{\ast}\left( {\mathbb F}_{\infty}\right) \to X$, 
the ${\mathcal C}$-invertibility of $p$ guarantees that there exists a surjective $\ast$-homomorphism
$h \colon Z_{C} \to X$ such that $g = h\circ p$.

In Section \ref{S:decomposition} we consider the case
${\mathcal C} = {\mathcal B}{\mathcal R}_{n}^{K}$, where $K > 0$ and ${\mathcal B}{\mathcal R}_{n}^{K}$
denotes the class of all separable unital $C^{\ast}$-algebras
of bounded rank with respect to $K$ at most $n$,
and prove that there indeed exists a ${\mathcal B}{\mathcal R}_{n}^{K}$-invertible
morphism
$p \colon C^{\ast}\left( {\mathbb F}_{\infty}\right) \to Z_{n}^{K}$
such that $Z_{n}^{K} \in {\mathcal B}{\mathcal R}_{n}^{K}$ which, as noted, implies that
$Z_{n}^{K}$ is a universal $C^{\ast}$-algebra in the class ${\mathcal B}{\mathcal R}_{n}^{K}$. It is interesting to note that not only $Z_{n}^{K}$ is universal in the above sense, but for every unital $\ast$-homomorphism $g \colon C^{\ast}(\mathbb{F}_{\infty}) \to \mathbb{C}$ the pushout $Z_{n}^{K}\bigstar_{\mathbb{C}}\mathbb{C}$, generated by $p$ and $g$, is also ${\mathcal B}{\mathcal R}_{n}^{K}$-universal. To see this take any separable unital $C^{\ast}$-algebra $X$ such that $\operatorname{br}_{K}(X) \leq n$ and consider the unital $\ast$-homomorphism $h \colon \mathbb{C} \to X$. Since $p$ is ${\mathcal B}{\mathcal R}_{n}^{K}$-invertible, there exists a unital $\ast$-homomorphism $\tilde{h} \colon Z_{n}^{K} \to X$ such that $\tilde{h}\circ p = h\circ g$. The homomorphisms $h$ and $\tilde{h}$ uniquely determine the unital $\ast$-homomorphim $\varphi \colon Z_{n}^{K}\bigstar_{\mathbb{C}}\mathbb{C} \to X$ as required.

%%%%%%%%%%%%%%%%%%%%%%%%%%%%%

\section{Preliminaries}\label{S:pre}
All $C^{\ast}$-algebras below are assumed to be unital.\
When we refer to a unital $C^{\ast}$-subalgebra of a
unital $C^{\ast}$-algebra we implicitly assume that the
inclusion is a unital $\ast$-homomorphism. The set of all
self-adjoint elements of a $C^{\ast}$-algebra $X$ is denoted
by $X_{sa}$. The product in the category of (unital) $C^{\ast}$-algebras, i.e. the $\ell^{\infty}$-direct sum, is denoted by $\prod\{ X_{t} \colon t \in T \}$.
For a given set $Y$ and a cardinal number $\tau$ the symbol
$\exp_{\tau}Y$ denotes the partially ordered (by inclusion) set of
all subsets of $Y$ of cardinality not exceeding $\tau$.

Recall that a direct system
${\mathcal S} = \{ X_{\alpha}, i_{\alpha}^{\beta}, A\}$ of
unital $C^{\ast}$-algebras consists of a partially ordered directed
indexing set $A$,
unital $C^{\ast}$-algebras $X_{\alpha}$, $\alpha \in A$, and 
unital $\ast$-homomorphisms
$i_{\alpha}^{\beta} \colon X_{\alpha} \to X_{\beta}$,
defined for each
pair of indexes $\alpha ,\beta \in A$ with $\alpha \leq \beta$,
and satisfying
the condition $i_{\alpha}^{\gamma} =
i_{\beta}^{\gamma}\circ i_{\alpha}^{\beta}$ for
each triple of indexes
$\alpha ,\beta ,\gamma \in A$ with $\alpha \leq \beta \leq \gamma$.
The (inductive) limit of the above direct
system is a unital $C^{\ast}$-algebra which is denoted by
$\varinjlim{\mathcal S}$. For each $\alpha \in A$ there
exists a unital $\ast$-homomorphism
$i_{\alpha} \colon X_{\alpha} \to \varinjlim{\mathcal S}$
which will be called
the $\alpha$-th limit homomorphism of $\mathcal S$. 

If $A^{\prime}$ is a
directed subset of the indexing set $A$, then the subsystem
$\{ X_{\alpha}, i_{\alpha}^{\beta}, A^{\prime}\}$ of
${\mathcal S}$ is denoted ${\mathcal S}|A^{\prime}$.

In Section \ref{S:decomposition} we use the concept of the direct
$C_{\tau}^{\ast}$-system introduced in \cite{chi991}.

\begin{defin}\label{D:smooth}
Let $\tau \geq \omega$
be a cardinal number. A direct system
${\mathcal S} = \{ X_{\alpha}, i_{\alpha}^{\beta}, A\}$ of
unital $C^{\ast}$-algebras $X_{\alpha}$
and unital $\ast$-ho\-mo\-morp\-hisms
$i_{\alpha}^{\beta} \colon X_{\alpha} \to X_{\beta}$ is called a {\em direct
$C_{\tau}^{\ast}$-system} if the
following conditions are satisfied:
\begin{itemize}
\item[(a)]
$A$ is a $\tau$-complete set (this means that for each chain $C$ of elements of the directed set $A$ with
${\mid}C{\mid} \leq \tau$, there exists an element
$\sup C$ in $A$; see \cite{book} for details). 
\item[(b)]
Density of $X_{\alpha}$ is at most $\tau$
(i.e. $d(X_{\alpha}) \leq \tau$), $\alpha \in A$.
\item[(c)]
The $\alpha$-th limit homomorphism
$i_{\alpha} \colon X_{\alpha} \to \varinjlim{\mathcal S}$ is
an injective $\ast$-ho\-mo\-mor\-phism for each $\alpha \in A$.
\item[(d)]
If $B = \{ \alpha_{t} \colon t \in T\} $ is a chain of
elements of $A$ with
$|T| \leq \tau$ and $\alpha = \sup B$, then the limit
homomorphism
$\varinjlim\{ i_{\alpha_{t}}^{\alpha} \colon t \in T\}
\colon \varinjlim\left({\mathcal S}|B\right)
\to X_{\alpha}$ is an isomorphism.
\end{itemize}
\end{defin}

\begin{pro}[Proposition 3.2, \cite{chi991}]\label{P:exists}
Let $\tau$ be an infinite cardinal number. Every unital
$C^{\ast}$-algebra $X$
can be represented as the limit of a direct
$C_{\tau}^{\ast}$-system
${\mathcal S}_{X} = \{ X_{\alpha}, i_{\alpha}^{\beta},
A \}$ where the indexing set $A$ coincides with $\exp_{\tau}Y$ for some (any)
dense subset $Y$ of $X$ with $|Y| = d(X)$.
\end{pro}

\begin{lem}[Lemma 3.3, \cite{chi991}]\label{L:strong}
If ${\mathcal S}_{X} = \{ X_{\alpha}, i_{\alpha}^{\beta}, A\}$
is a direct $C_{\tau}^{\ast}$-system, then 
\[ \varinjlim{\mathcal S}_{X} = \bigcup\{
i_{\alpha}(X_{\alpha}) \colon \alpha \in A\} .\]
\end{lem}

%%%%%%%%%%%%%%%%%%%%%%%%%%%%%%%%%%%%%%%%%%%%%%%%%%%%%%%%%%%
%%%%%%%%%%%%%%%%%%%%%%%%%%%%%%%%%%%%%%%%%%%%%%%%%%%%%%%%

The real rank of a unital $C^{\ast}$-algebra $X$, denoted by $\operatorname{rr}(X)$,
is defined as follows \cite{brownped91}. We say that $\operatorname{rr}(X) \leq n$ if for each
$(n+1)$-tuple $(x_{1},\dots ,x_{n+1})$ of self-adjoint elements
in $X$ and every $\epsilon > 0$,
there exists an $(n+1)$-tuple $(y_{1},\dots ,y_{n+1})$ in $X_{sa}$
such that $\sum_{k=1}^{n+1} y_{k}^{2}$ is invertible and
$\left\| \sum_{k=1}^{n+1} (x_{k}-y_{k})^{2}\right\| < \epsilon$ .

\begin{lem}\label{L:definition}
The following conditions are equivalent for a unital $C^{\ast}$-algebra $X$;
\begin{itemize}
\item[(i)]
$\operatorname{rr}(X) \leq n$;
\item[(ii)]
for each $(n+1)$-tuple $(x_{1},\dots ,x_{n+1})$ 
in $X_{sa}$ and for each $\epsilon > 0$,
there exists an $(n+1)$-tuple $(y_{1},\dots ,y_{n+1})$ in $X_{sa}$
such that $\sum_{k=1}^{n+1} y_{k}^{2}$ is invertible and
$\left\| x_{k}-y_{k}\right\| < \epsilon$
for each $k = 1,2,\dots ,n+1$.
\item[(iii)]
for each $(n+1)$-tuple $(x_{1},\dots ,x_{n+1})$ 
in $X_{sa}$, with $\| x_{k}\| =1$ for each $k = 1,\dots , n+1$, and for each $\epsilon > 0$,
there exists an $(n+1)$-tuple $(y_{1},\dots ,y_{n+1})$ in $X_{sa}$
such that $\sum_{k=1}^{n+1} y_{k}^{2}$ is invertible and
$\left\| x_{k}-y_{k}\right\| < \epsilon$
for each $k = 1,\dots ,n+1$.
\end{itemize}
\end{lem}
\begin{proof}
(i)$\Longrightarrow$(ii). Let  $(x_{1},\dots ,x_{n+1})$  be an
$(n+1)$-tuple in $X_{sa}$ and $\epsilon > 0$. By (i), there exists an $(n+1)$-tuple $(y_{1},\dots ,y_{n+1})$ in $X_{sa}$
such that $\sum_{k=1}^{n+1} y_{k}^{2}$ is invertible and
$\left\| \sum_{k=1}^{n+1} (x_{k}-y_{k})^{2}\right\| < \epsilon^{2}$.

Since $x_{k} -y_{k} \in X_{sa}$, it follows (\cite[2.2.4 Theorem]{murphy}) that
$(x_{k}-y_{k})^{2} \geq 0$ for each $k = 1,\dots ,n+1$. Then, by \cite[2.2.3 Lemma]{murphy},
$\sum_{k=1}^{n+1} (x_{k}-y_{k})^{2} \geq 0$. Note also that 
$\sum_{i=1}^{n+1} (x_{i}-y_{i})^{2} - (x_{k} -y_{k})^{2} = \sum_{i=1,i \neq k}^{n+1} (x_{i}-y_{i})^{2} \geq 0, k = 1,\dots ,n+1,$ which guarantees that $(x_{k} -y_{k})^{2} \leq \sum_{i=1}^{n+1} (x_{i}-y_{i})^{2}$ for each $k = 1,\dots ,n+1$. By \cite[2.2.5 Theorem]{murphy}, 
$\| x_{k}-y_{k}\|^{2}  = \| (x_{k}-y_{k})^{2}\| \leq \left\| \sum_{i=1}^{n+1} (x_{i}-y_{i})^{2}\right\| < \epsilon^{2}$. Consequently $\| x_{k}-y_{k}\| < \epsilon$, $k = 1,\dots ,n+1$. This shows that condition (ii) is satisfied.

The implication (ii)$\Longrightarrow$(iii) is trivial.

(iii)$\Longrightarrow$(i). Let $(x_{1},\dots ,x_{n+1})$ be an
$(n+1)$-tuple in $X_{sa}$ and $\epsilon > 0$. We can assume that each $x_{i}\neq 0$. Consider the $(n+1)$-tuple $\left( \frac{x_{1}}{\| x_{1}\|},\dots ,\frac{x_{n+1}}{\| x_{n+1}\|}\right)$. By (iii), there exists an $(n+1)$-tuple $( \bar{y}_{1},\dots ,\bar{y}_{n+1})$ 
in $X_{sa}$ such that $\sum_{k=1}^{n+1} \bar{y}_{k}^{2}$ is invertible and
$\left\| \frac{x_{k}}{\| x_{k}\|}-\bar{y}_{k}\right\| <  \sqrt{\frac{M\cdot \epsilon}{n+1}}$, $k = 1,2,\dots ,n+1$, where $M = \frac{1}{\max\{ \|x_{k}\|^{2} \colon k = 1,\dots ,n+1 \}}$. Now let $y_{k} = \| x_{k}\| \cdot \bar{y}_{k}$ for each $k = 1,\dots ,n+1$. Then we have $\| (x_{k} -y_{k})^{2}\| = \| (x_{k}-y_{k})(x_{k}-y_{k})^{\ast}\| = \| x_{k}-y_{k}\|^{2} = \left\| \| x_{k}\|\cdot \frac{x_{k}}{\| x_{k}\|} -\| x_{k}\|\cdot\bar{y}_{k}\right\|^{2} = 
\|x_{k}\|^{2}\left\| \frac{x_{k}}{\| x_{k}\|} - \bar{y}_{k}\right\|^{2}< \frac{\epsilon}{n+1}, k = 1,\dots ,n+1$. Consquently $\left\| \sum_{k=1}^{n+1}(x_{k}-y_{k})^{2}\right\| \leq \sum_{k=1}^{n+1}\| (x_{k}-y_{k})^{2}\| < \epsilon$, which shows that $\operatorname{rr}(X) \leq n$.
\end{proof}

%%%%%%%%%%%%%%%%%%%%%%%%%%%%%%%%%%%%%%%%%%%%%%%%%%%%%%%%%%%%
%%%%%%%%%%%%%%%%%%%%%%%%%%%%%%%%%%%%%%%%%%%%%%%%%%%%%%%%%%%

\section{Bounded rank, real rank and joint spectra}\label{S:rjs}
One of our main goals in this Section is to define the bounded rank of unital $C^{\ast}$-algebras. Formal definition makes it not very easy to see how this new concept differs from the concept of the real rank. Intuitive geometrical interpretation does not really help. Moreover, we show that in the presence of nice joint spectrum real and bounded ranks are indeed identical. Such situations include the commutative case as well as the case of real rank zero.

\subsection{Axioms of generalized joint spectra}\label{SS:axioms}

Let $\mathcal{K}(\mathbb{C}^{m})$ denote the collection of compact subsets of $\mathbb{C}^{m}$. Let $M \in \N \cup \{\infty\}$. A generalized joint $M$-spectrum (or simply spectrum) on an unital $C^{\ast}$-algebra $X$ is a collection $\{\widetilde{\sigma}_{m} \colon m \leq M\}$  of maps $\widetilde{\sigma}_{m} \colon  X^{m} \to \mathcal{K}(\mathbb{C}^{m})$
satisfying conditions (I)--(III) below:

\begin{itemize}
\item[(I)]
$\displaystyle \widetilde{\sigma}_{m}(x_{1},\dots ,x_{m})$ is a non-empty compact subset of $\mathbb{C}^{m}$.
\item[(II)]
If $x \in X$, then $\widetilde{\sigma}_{1}(x) = \sigma (x)$, where $\sigma (x)$ denotes the usual spectrum of the element $x$.
\item[(III)]
$\displaystyle p(\widetilde{\sigma}_{m}(x_{1},\dots ,x_{m})) = \widetilde{\sigma}_{n}(p(x_{1},\dots ,x_{m}))$, where $p = (p_{1},\dots ,p_{n}) \colon \mathbb{C}^{m} \to\mathbb{C}^{n}$ is a polynomial mapping.
\end{itemize}
\noindent We also consider the following two properties:
\begin{itemize}
\item[(IV)]
There exists a constant $K >0$ (called a polynomial spectral constant) such that for any rational $\delta >0$, for any $m$-tuple $(x_{1},\dots ,x_{m})$ in $X$ (resp., in $X_{sa}$) and for any polynomial mapping $p = (p_{1},\dots ,p_{m}) \colon \mathbb{C}^{m} \to \mathbb{C}^{m}$ (resp., having real coefficients) with $\| p - \operatorname{id}_{\mathbb{C}^{m}}\| \leq K\cdot\delta$, there exists an $m$-tuple $(y_{1},\dots ,y_{m})$ in $X$ (resp., in $X_{sa}$) such that $\widetilde{\sigma}(y_{1},\dots ,y_{m}) = p(\widetilde{\sigma}_{m}(x_{1},\dots ,x_{m}))$ and $\| x_{k}-y_{k}\| \leq \delta$ for each $k=1,\dots ,m$.
\item[(V)]
There exists a constant $K > 0$ (called a general spectral constant) such that for any rational $\delta > 0$, for any $m$-tuple $(x_{1},\dots ,x_{m})$ in $X$ (resp., in $X_{sa}$) and for any map $f $ from $ \widetilde{\sigma}_{m}(x_{1},\dots ,x_{m})$ into $ \mathbb{C}^{m}$ (resp., into $\mathbb{R}^{m} \subseteq \mathbb{C}^{m}$) with $\| f - \operatorname{id}_{\widetilde{\sigma}_{m}(x_{1},\dots ,x_{m}) }\| \leq K\cdot\delta$, there exists an $m$-tuple $(y_{1},\dots ,y_{m})$ in $X$ (resp., in $X_{sa}$) such that $\widetilde{\sigma}(y_{1},\dots ,y_{m}) = f(\widetilde{\sigma}_{m}(x_{1},\dots ,x_{m}))$ and $\| x_{k}-y_{k}\| \leq \delta$ for each $k=1,\dots ,m$.
\end{itemize}

Few comments are in order: 

\begin{itemize}
\item[(A)]
Concerning property (I), it is not usually requested (in a much general setting though)  the spectrum of every tuple of non-commuting elements to be non-empty. 
\item[(B)]
Properties (I)--(III) are standard axioms \cite{zelazko} of joint spectra (for commuting tuples) in Banach algebras. Property (III) is known as the spectral mapping property of $\widetilde{\sigma}$.  
\item[(C)]
It is easy to see \cite{zelazko} that properties (II) and (III) imply the inclusion $\widetilde{\sigma}_{m}(x_{1},\dots ,x_{m}) \subseteq \prod\{ \sigma (x_{k}) \colon k = 1,\dots ,m\}$ for any $m$-tuple $(x_{1},\dots ,x_{m})$ in $X$.
\item[(D)]
Note that if, in property (III), the $m$-tuple $(x_{1},\dots ,x_{m})$ consists of self-adjoint elements and the polynomials $p_{k}$, $k=1,\dots ,m$, have real coefficients, then the $m$-tuple $p(x_{1},\dots ,x_{m}) = (p_{1}(x_{1},\dots ,x_{m}),\dots ,p_{m}(x_{1},\dots ,x_{m}))$ also consists of self-adjoint elements.
\item[(E)]
In property (IV), we do not specify whether the $m$-tuple $(y_{1},\dots ,y_{m})$ is obtained as the image of the $m$-tuple 
$(x_{1},\dots ,x_{m})$
under the polynomial $p$. Nevertheless, in light of (D), we require that all $y_{k}$'s are self-adjoint provided that all $x_{k}$'s are self-adjoint and all the coefficients of the polynomials $p_{k}$ are real. 
\item[(F)]
Similar comment with respect to property (V). If all $x_{k}$'s are self-adjoint and $f(\widetilde{\sigma}_{m}(x_{1},\dots ,x_{m})) \subseteq \mathbb{R}^{m} \subseteq \mathbb{C}^{m}$, then all $y_{k}$'s are also assumed to be self-adjoint.
\item[(G)]
If $T$ is a compact space and $(x_{1},\dots ,x_{m})$ is an $m$-tuple in $C(T)$, let define $\widetilde{\sigma}_{m}(x_{1},\dots ,x_{m}) = \triangle\{ x_{k} \colon k = 1,\dots ,m\}(T)$, where $\triangle\{ x_{k} \colon k = 1,\dots ,m\} (t) = (x_{1}(t),\dots ,x_{m}(t))$ for each $t \in T$. So obtained generalized joint spectrum has all the properties (I)--(V) with $K=1$. 
\item[(H)]
For any  $C^{\ast}$-algebra $X$ the usual
spectrum $\sigma (x)$, $x\in X$,  provides an example of a generalized joint $1$-spectrum. This follows from spectral mapping theorem and considerations related to the functional calculus.
\end{itemize}

\vspace{0.4in}
%%%%%%%%%%%%%%%%%%%%%%%%%%%%%%%%%%%%%%%%%%%%%%%%%%%%
%%%%%%%%%%%%%%%%%%%%%%%%%%%%%%%%%%%%%%%%%%%%%%%%%%%%%

\subsection{Real rank in terms of generalized joint spectra}\label{SS:rr}
Before introducing the concept of the bounded rank we illustrate how the usual real rank can be characterized in terms of generalized joint spectra. This result will be used later to compare real and bounded ranks.
\begin{lem}\label{L:rjs}
Let $X$ be a unital $C^{\ast}$-algebra with a generalized joint $m$-spectrum. Then the following conditions are equivalent for any $m$-tuple $(y_{1},\dots ,y_{m})$ of self-adjoint elements in $X$:
\begin{itemize}
\item[(a)]
The element $\sum_{k=1}^{m}y_{k}^{2}$ is invertible.
\item[(b)]
$\mathbf{0} \not\in \widetilde{\sigma}_{m}(y_{1},\dots ,y_{m})$.
\end{itemize}
\end{lem}
\begin{proof}
 Since $y_{k} = y_{k}^{\ast}$ it follows that $\widetilde{\sigma}_{1}(y_{k}) = \sigma (y_{k}) \subseteq \mathbb{R}$ and consequently, by properties (II) and (III) (see comment (C)) of the joint spectrum, $\widetilde{\sigma}_{m}(y_{1},\dots ,y_{m}) \subseteq \mathbb{R}^{m}$. Consider the polynomial $p(t_{1},\dots ,t_{m}) = \sum_{k=1}^{m}t_{k}^{2}$. By property (III),\\ $p\left(\widetilde{\sigma}_{m}(y_{1},\dots ,y_{m})\right) =  \widetilde{\sigma}_{1}(p(y_{1},\dots ,y_{m})) =  \sigma\left( \sum_{k=1}^{m}y_{k}^{2}\right) \subseteq [0,\infty )$. This implies that  $\mathbf{0} \not\in \widetilde{\sigma}_{m}(y_{1},\dots ,y_{m})$ if and only if
$\sigma\left( \sum_{k=1}^{m}y_{k}^{2}\right) \subseteq (0,\infty ) $. Finally, note that invertibility of  $\sum_{k=1}^{m}y_{k}^{2}$ is equivalent to $\sigma \left(\sum_{k=1}^{m}y_{k}^{2}\right) \subseteq (0,\infty)$.
\end{proof}

As a corollary we obtain the following statement

\begin{pro}\label{P:rjs}
If there exists a generalized joint $(n+1)$-spectrum on a unital $C^{\ast}$-algebra $X$, then the following conditions are equivalent:
\begin{itemize}
\item[(a)]
$\operatorname{rr}(X) \leq n$;
\item[(b)]
For every very $(n+1)$-tuple $(x_{1}, \dots ,x_{n+1})$ of self-adjoint elements in $X$ and for every $\epsilon > 0$ there exists an $(n+1)$-tuple $(y_{1}, \dots ,y_{n+1})$ of self-adjoint elements such that
\begin{itemize}
\item[(i)]
$\| x_{k} -y_{k}\| < \epsilon$ for each $k = 1,\dots ,n+1$;
\item[(ii)]
$\mathbf{0} \not\in \widetilde{\sigma}_{n+1}(y_{1},\dots ,y_{n+1})$.
\end{itemize}
\end{itemize}
\end{pro}
%%%%%%%%%%%%%%%%%%%%%%%%%%%%%%%%%%%%%%%%%%%%%%%%
%%%%%%%%%%%%%%%%%%%%%%%%%%%%%%%%%%%%%%%%%%%%%%%

\vspace{0.4in}

\subsection{Bounded rank -- main definitions}\label{SS:br}

We begin this section with our main definitions.

\begin{defin}\label{D:modifiedunessential}
Let $K >0$. We say that an $m$-tuple $( y_{1},\dots ,y_{m})$ of self-adjoint elements of a
unital $C^{\ast}$-algebra $X$ is $K$-{\em unessential} if for every rational $\delta > 0$ there exists
an $m$-tuple $(z_{1},\dots ,z_{m})$ of self-adjoint elements of $X$ satisfying the following conditions:
\begin{itemize}
\item[(a)]
$\| y_{k}-z_{k}\| \leq\delta$ for each $k = 1,\dots ,m$,
\item[(b)]
The element $\sum_{k=1}^{m}z_{k}^{2}$ is invertible and $\left\| \left(\sum_{k=1}^{m}z_{k}^{2}\right)^{-1}\right\| \leq \frac{1}{K\cdot \delta^{2}}$.
\end{itemize}
$1$-unessential tuples are referred as {\em unessential}.
\end{defin}

\begin{rem}\label{R:different}
Obviously if $K_{1} \leq K_{2}$, then every $K_{2}$-unessential $m$-tuple is $K_{1}$-unessential.
\end{rem}

\begin{defin}\label{D:modifiedbrr}
Let $K>0$. We say that the {\em bounded rank} of a unital $C^{\ast}$-algebra $X$ with respect to $K$ 
does not exceed $n$ (notation: $\operatorname{br}_{K}(X) \leq n$) if for any 
$(n+1)$-tuple $(x_{1},\dots ,x_{n+1})$ of self-adjoint elements of $X$ and for any $\epsilon > 0$ there exists a $K$-unessential $(n+1)$-tuple $(y_{1},\dots ,y_{n+1})$ in $X$ such that $\| x_{k}-y_{k}\| < \epsilon$ for each $k = 1,\dots ,n+1$. For simplicity $\operatorname{br}_{1}(X)$  is denoted by $\operatorname{br}(X)$ and it is called a {\em bounded rank}.
\end{defin}

We record the following statement for future references.

\begin{pro}\label{P:comparison1}
$\operatorname{rr}(X) \leq \operatorname{br}_{K}(X)$ for any unital $C^{\ast}$-algebra $X$ and for any $K>0$.
\end{pro}
\begin{proof}
Let $\operatorname{br}_{K}(X) = n$ and  $(x_{1},\dots ,x_{n+1})$ be an $(n+1)$-tuple of self-adjoint elements in $X$. Let also $\epsilon > 0$.  
Since $\operatorname{br}_{K}(X) = n$, there exists a $K$-unessential $(n+1)$-tuple $(y_{1},\dots ,y_{n+1})$ such that $\| x_{k}-y_{k}\| < \frac{\epsilon}{2}$ for each $k=1,\dots ,n+1$. This in turn means that for a rational number $\delta$ with $\delta \leq \frac{\epsilon}{2}$ there is an $(n+1)$-tuple $(z_{1},\dots ,z_{n+1})$ of self-adjoint elements such that
\begin{itemize}
\item[(a)]
$\|y_{k}-z_{k}\| \leq \delta$ for each $k = 1,\dots ,n+1$,
\item[(b)]
The element $\sum_{k=1}^{n+1}z_{k}^{2}$ is invertible and $\left\| \left(\sum_{k=1}^{n+1}z_{k}^{2}\right)^{-1}\right\| \leq \frac{1}{K\cdot \delta^2}$.
\end{itemize}
Clearly $\| x_{k}-z_{k}\| \leq \| x_{k}-y_{k}\| +\|y_{k}-z_{k}\| < \frac{\epsilon}{2}+\delta \leq\epsilon$, $k=1,\dots ,n+1$. According to (b), $\sum_{k=1}^{n+1}z_{k}^{2}$ is invertible which shows that $\operatorname{rr}(X) \leq n$.
\end{proof}

%%%%%%%%%%%%%%%%%%%%%%%%%%%%%%%%%%%%%%%%
%%%%%%%%%%%%%%%%%%%%%%%%%%%%%%%%%%%%%%%%%
\vspace{0.4in}

\subsubsection{Bounded rank and polynomial spectral constant of a joint spectrum}\label{SSS:pspectral}
In this subsection we investigate relations between bounded and real ranks in the presence of a generalized joint spectrum satisfying condition (IV).

\begin{lem}\label{L:bremma}
Let $X$ be a unital $C^{\ast}$-algebra with a generalized joint $m$-spectrum satisfying condition $(IV)$. If $K$ is a polynomial spectral constant of the joint spectrum, then for any $m$-tuple $(y_{1},\dots ,y_{m})$ of self-adjoint elements of $X$ we have {\em (a)}$\Longrightarrow${\em (b)}$\Longrightarrow${\em (c)}, where
\begin{itemize}
\item[(a)]
$\mathbf{0} \not\in \widetilde{\sigma}_{m}(y_{1},\dots ,y_{m})$.
\item[(b)]
For any rational $\delta >0$ there exists an $m$-tuple $(z_{1},\dots z_{m})$ of self-adjoint elements in $X$ such that
\begin{itemize}
\item[(i)]
$\| z_{k}-y_{k}\| \leq \delta$ for each $k = 1,\dots ,m$.
\item[(ii)]
$\widetilde{\sigma}_{m}(z_{1},\dots ,z_{m})) \subseteq \mathbb{R}^{m} \setminus O\left(\mathbf{0},\frac{K}{2}\cdot \delta\right)$.
\end{itemize}
\item[(c)]
$(y_{1},\dots ,y_{m})$ is $\frac{K^{2}}{4}$-unessential. 
\end{itemize}
\end{lem}
\begin{proof}
(a)$\Longrightarrow$(b).  Choose a rational number $\delta > 0$. Our goal is to find an $m$-tuple $(z_{1},\dots ,z_{m})$ such that
\begin{itemize}
\item[(i)]
$\| y_{k}-z_{k}\| \leq \delta$ for each $k=1,\dots ,m$.
\item[(ii)]
$\widetilde{\sigma}_{m}(z_{1},\dots ,z_{m})) \subseteq \mathbb{R}^{m} \setminus O\left(\mathbf{0},\frac{K}{2}\cdot \delta\right)$. 
\end{itemize}

 Since $\mathbf{0} \not\in \widetilde{\sigma}_{m}(y_{1},\dots ,y_{m})$, there exists a number $C > 0$ such that $O(\mathbf{0},C) \cap \widetilde{\sigma}_{m}(y_{1},\dots ,y_{m}) = \emptyset$. Without loss of generality we may assume that $C \leq \frac{K\cdot \delta}{2}$.

Consider the function $r_{\delta} \colon \mathbb{R} \to \mathbb{R}$ defined by letting
\[
r_{\delta}(t) = 
\begin{cases}
\displaystyle t,\;\text{if}\; |t| \geq \frac{K\cdot \delta}{2}+C  ,\\
\displaystyle \frac{K\cdot \delta}{2}+C ,\; \text{if}\; C \leq t \leq \frac{K\cdot \delta}{2}+C ,\\
\displaystyle -\left(\frac{K\cdot \delta}{2}+C\right) ,\; \text{if}\;  -\frac{K\cdot \delta}{2}-C \leq t \leq -C,\\
\displaystyle \left(\frac{K\cdot \delta}{2C}+1\right)\cdot t ,\; \text{if}\; |t| \leq C .\\
\end{cases}
\]

\noindent Let $p_{\delta} \colon \mathbb{R} \to \mathbb{R}$ be a polynomial such that  $| p_{\delta}-r_{\delta}| \leq C$. Next consider the polynomial $\tilde{p}_{\delta}(t_{1},\dots ,t_{m}) = (p_{\delta}(t_{1}),\dots ,p_{\delta}(t_{m})) \colon \mathbb{R}^{m} \to \mathbb{R}^{m}$. 
Note that 
\begin{equation}\label{E:close}
\| \tilde{p}_{\delta}-\operatorname{id}\| = |p_{\delta} - \operatorname{id}| \leq
 |p_{\delta}-r_{\delta}| +| r_{\delta}- \operatorname{id}| \leq C +\frac{K\cdot\delta}{2} \leq K\cdot \delta .
\end{equation}

By property (IV), there exists an $m$-tuple $(z_{1},\dots ,z_{k})$ such that $\| z_{k} -y_{k}\| \leq \delta$ for each $k = 1,\dots ,m$ and 
\begin{equation}\label{E:contained}
\widetilde{\sigma}_{m}(z_{1},\dots ,z_{m}) = \tilde{p}_{\delta}(\widetilde{\sigma}_{m}(y_{1},\dots ,y_{m})) .
\end{equation}

If $(t_{1},\dots ,t_{m}) \in \widetilde{\sigma}_{m}(y_{1},\dots ,y_{m})$, then for at least one $k = 1,\dots ,m$, we must have $|t_{k}| \geq C$ (recall that $O(\mathbf{0},C) \cap \widetilde{\sigma}_{m}(y_{1},\dots ,y_{m}) = \emptyset$). Hence, $|r_{\delta}(t_{k})| \geq \frac{K\cdot \delta}{2}+C$. Consequently there are two possibilities
\smallskip

\[
\begin{cases}
\displaystyle r_{\delta}(t_{k}) -C \leq p_{\delta}(t_{k}) \leq r_{\delta}(t_{k}) +C\\
\displaystyle r_{\delta}(t_{k}) \geq \frac{K\cdot \delta}{2}+C\\
\end{cases}
\Longrightarrow 
 p_{\delta}(t_{k}) \geq r_{\delta}(t_{k}) -C \geq  \frac{K\cdot \delta}{2}
\]
\smallskip

\noindent or

\[
\begin{cases}
\displaystyle r_{\delta}(t_{k}) -C \leq p_{\delta}(t_{k}) \leq r_{\delta}(t_{k}) +C\\
\displaystyle r_{\delta}(t_{k}) \leq -\frac{K\cdot \delta}{2}-C\\
\end{cases} 
\Longrightarrow p_{\delta}(t_{k}) \leq r_{\delta}(t_{k}) +C \leq -\frac{K\cdot \delta}{2}
\] 
\smallskip

\noindent which simply means that $|p_{\delta}(t_{k})| \geq \frac{K\cdot \delta}{2}$.  Therefore,
\begin{equation}\label{E:below}
 \| \tilde{p}_{\delta}(t_{1},\dots ,t_{m})\| = \max\{ |p_{\delta}(t_{i})| \colon i = 1,\dots ,m\} \geq |p_{\delta}(t_{k})| \geq \frac{K\cdot \delta}{2}.
\end{equation}

By (\ref{E:contained}) and (\ref{E:below}), 

\begin{equation}\label{E:outside}
\widetilde{\sigma}_{m}(z_{1},\dots ,z_{m}) \subseteq \mathbb{R}^{m} \setminus O\left(\mathbf{0},\frac{K\cdot \delta}{2}\right)
\end{equation}

(b)$\Longrightarrow$(c).  Let $\delta>0$ be rational. We fix an $m$-tuple $(z_{1},\dots ,z_{m})$ of self-adjoint elements of $X$ satisfying conditions (i) and (ii) from (b). In order to prove that $(y_{1},\dots ,y_{m})$ is $\frac{K^{2}}{4}$-unessential, we need only to check that  $\left\| \left(\sum_{k=1}^{m}z_{k}^{2}\right)^{-1}\right\| \leq\frac{4}{K^{2}\cdot\delta^{2}}$. To this end, consider the polynomial $\lambda (t_{1},\dots ,t_{m}) = \sum_{k=1}^{m}t_{k}^{2}$. According to Property (III), we have 
 \begin{equation}\label{E:contained1}
\sigma\left( \sum_{k=1}^{m}z_{k}^{2} \right) = \lambda \left(\widetilde{\sigma}_{m}(z_{1},\dots ,z_{m})\right) \subseteq \lambda\left(\mathbb{R}^{m} \setminus O\left(\mathbf{0},\frac{K\cdot \delta}{2}\right)\right) ,
\end{equation}

\noindent  If $(t_{1},\dots ,t_{m}) \in \mathbb{R}^{m} \setminus O\left(\mathbf{0},\frac{K\cdot \delta}{2}\right)$, then there exists $i = 1,\dots ,m$ such that $|t_{i}| \geq \frac{K\cdot \delta}{2}$ and consequently,
$\sum_{k=1}^{m}t_{k}^{2} \geq \frac{K^{2}}{4}\cdot \delta^{2}$. By (\ref{E:contained1}), we finally obtain
\[\sigma\left( \sum_{k=1}^{m}z_{k}^{2} \right) \subseteq \left[ \frac{K^{2}}{4}\cdot \delta^{2},\infty\right) \Longrightarrow \left\| \left(\sum_{k=1}^{m}z_{k}^{2}\right)^{-1}\right\| \leq
\frac{4}{K^{2}\cdot\delta^{2}}.\]
\end{proof}

As a corollary of Lemma \ref{L:bremma} we obtain the following statement.

\begin{pro}\label{P:maincomparison}
If there exists a generalized joint $(n+1)$-spectrum on a unital $C^{\ast}$-algebra $X$ satisfying conditions $(IV)$, then the following  are equivalent:
\begin{itemize}
\item[(a)]
$\operatorname{rr}(X) \leq n$;
\item[(b)]
$\operatorname{br}_{\frac{K^2}{4}}(X) \leq n$, where $K$ is a polynomial spectral constant of the joint spectrum. 
\end{itemize}
\end{pro}
\begin{proof}
The implication (b)$\Longrightarrow$(a) follows from Proposition \ref{P:comparison1}. To prove the implication (a)$\Longrightarrow$(b), let 
 $\epsilon > 0$ and  $(x_{1},\dots ,x_{n+1})$ be an $(n+1)$-tuple of self-adjoint elements of $X$. By Proposition \ref{P:rjs}, there exists an $(n+1)$-tuple $(y_{1},\dots ,y_{n+1})$ of self-adjoint elements of $X$ such that $\| x_{k}-y_{k}\| < \epsilon$ for each $k = 1,\dots ,n+1$ and $\mathbf{0} \not\in \widetilde{\sigma}_{n+1}(y_{1},\dots ,y_{n+1})$. Then Lemma
\ref{L:bremma} yields that the $(n+1)$-tuple $(y_{1},\dots ,y_{n+1})$ is $\frac{K^2}{4}$-unessential. Hence, $\operatorname{br}_{\frac{K^2}{4}}(X) \leq n$.
\end{proof}

%%%%%%%%%%%%%%%%%%%%%%%%%%%%%%%%%%%%%%%%
%%%%%%%%%%%%%%%%%%%%%%%%%%%%%%%%%%%%%%%%%

\vspace{0.4in}

\subsubsection{Bounded rank and a general spectral constant of a joint spectrum}\label{SSS:spectral}
In this subsection we investigate relations between bounded and real ranks in the presence of a generalized joint spectrum satisfying condition (V).

\begin{lem}\label{L:bremma2}
Let $X$ be a unital  $C^{\ast}$-algebra with a generalized joint $m$-spectrum satisfying condition {\em (V)}. If $K$ is a general spectral constant of the joint spectrum, then for any $m$-tuple $(y_{1},\dots ,y_{m})$ of self-adjoint elements of $X$ we have {\em (a)}$\Longrightarrow${\em (b)}$\Longrightarrow${\em (c)}:
\begin{itemize}
\item[(a)]
$\mathbf{0} \not\in \widetilde{\sigma}_{m}(y_{1},\dots ,y_{m})$.
\item[(b)]
For any rational $\delta >0$ there exists an $m$-tuple $(z_{1},\dots z_{m})$ in $X$ such that
\begin{itemize}
\item[(i)]
$\| z_{k}-y_{k}\| \leq \delta$ for each $k = 1,\dots ,m$.
\item[(ii)]
$\widetilde{\sigma}_{m}(z_{1},\dots ,z_{m})) \subseteq \mathbb{R}^{m} \setminus O\left(\mathbf{0},K\cdot \delta\right)$.
\end{itemize}
\item[(c)]
$(y_{1},\dots ,y_{m})$ is $\displaystyle K^{2}$-unessential.
\end{itemize}
\end{lem}
\begin{proof}
(a)$\Longrightarrow$(b). For a given rational $\delta > 0$ we need to find an $m$-tuple $(z_{1},\dots ,z_{m})$ such that
\begin{itemize}
\item[(i)]
$\| y_{k}-z_{k}\| \leq \delta$ for each $k=1,\dots ,m$.
\item[(ii)]
$\widetilde{\sigma}_{m}(z_{1},\dots ,z_{m})) \subseteq \mathbb{R}^{m} \setminus O\left(\mathbf{0},K\cdot \delta\right)$.
\end{itemize}

\noindent Consider the function $r_{\delta} \colon \mathbb{R}^{m}\setminus \{ \mathbf{0}\} \to \mathbb{R}^{m}$, defined by letting
\[
r_{\delta}(\mathbf{t}) = 
\begin{cases}
\mathbf{t},\;\text{if}\; \| \mathbf{t}\| \geq K\cdot \delta  ,\\
\frac{\mathbf{t}K\cdot \delta}{\| \mathbf{t}\|} ,\; \text{if}\; 0 < \| \mathbf{t}\| \leq K\cdot \delta ,\\
\end{cases}
\]

\noindent where $\mathbf{t} = (t_{1},\dots ,t_{m}) \in \mathbb{R}^{m}\setminus \{ \mathbf{0}\}$. Since $\mathbf{0} \not\in \widetilde{\sigma}_{m}(y_{1},\dots ,y_{m})$, the restriction map $r_{\delta}|\widetilde{\sigma}_{m}(y_{1},\dots ,y_{m}) \colon \widetilde{\sigma}_{m}(y_{1},\dots ,y_{m}) \to \mathbb{R}^{m}$ is well defined. Note that $\displaystyle r_{\delta}(\widetilde{\sigma}_{m}(y_{1},\dots ,y_{m})) \subseteq \mathbb{R}^{m}\setminus O(\mathbf{0},K\cdot \delta )$. Since $\| r_{\delta}|\widetilde{\sigma}_{m}(y_{1},\dots ,y_{m}) -\operatorname{id}_{\widetilde{\sigma}_{m}(y_{1},\dots ,y_{m})}\| \leq K\cdot \delta$,  by property (V), there exists an $m$-tuple $(z_{1},\dots ,z_{m})$ in $X$ such that $\| z_{k} -y_{k}\| \leq \delta$ for each $k = 1,\dots ,m$ and $\widetilde{\sigma}_{m}(z_{1},\dots ,z_{m}) = r_{\delta}(\widetilde{\sigma}_{m}(y_{1},\dots ,y_{m}))$. It only remains to note that
\[ \widetilde{\sigma}_{m}(z_{1},\dots ,z_{m}) = r_{\delta}(\widetilde{\sigma}_{m}(y_{1},\dots ,y_{m})) \subseteq \mathbb{R}^{m} \setminus O\left( \mathbf{0},K\cdot \delta \right) .\]

(b)$\Longrightarrow$(c). For a given rational $\delta>0$ we fix an $m$-tuple $(z_{1},\dots z_{m})$ of self-adjoint elements of $X$ satisfying conditions (i) and (ii) of (b).  
As in the proof of Lemma \ref{L:bremma}, implication (b)$\Longrightarrow$(c), 
\begin{equation}\label{E:contained2}
\sigma\left( \sum_{k=1}^{m}z_{k}^{2} \right) = \lambda \left(\widetilde{\sigma}_{m}(z_{1},\dots ,z_{m})\right) \subseteq \lambda\left(\mathbb{R}^{m} \setminus O\left(\mathbf{0},K\cdot \delta\right)\right) ,
\end{equation}

\noindent where $\lambda (t_{1},\dots ,t_{m}) = \sum_{k=1}^{m}t_{k}^{2}$. Note that if $(t_{1},\dots ,t_{m}) \in \mathbb{R}^{m} \setminus O\left(\mathbf{0},K\cdot \delta\right)$, then there exists $i = 1,\dots ,m$ such that $|t_{i}| \geq K\cdot \delta$ and consequently,
$\sum_{k=1}^{m}t_{k}^{2} \geq K^{2}\cdot \delta^{2}$. By (\ref{E:contained2}), 
\[\sigma\left( \sum_{k=1}^{m}z_{k}^{2} \right) \subseteq \lambda\left(\mathbb{R}^{m} \setminus O\left(\mathbf{0},K\cdot \delta\right)\right) \subseteq \left[ K^{2}\cdot \delta^{2},\infty\right) \Longrightarrow \left\| \left(\sum_{k=1}^{m}z_{k}^{2}\right)^{-1}\right\| \leq \frac{1}{K^{2}\cdot \delta^{2}}.\]
\end{proof}

As a corollary of Lemma \ref{L:bremma2} we obtain the following statement.

\begin{pro}\label{P:maincomparison2}
If there exists a generalized joint $(n+1)$-spectrum on the unital $C^{\ast}$-algebra $X$ satisfying condition  {\em (V)}, then the following  are equivalent:
\begin{itemize}
\item[(a)]
$\operatorname{rr}(X) \leq n$;
\item[(b)]
$\displaystyle \operatorname{br}_{K^2}(X) \leq n$, where $K$ is a general spectral constant of the joint spectrum. 
\end{itemize}
\end{pro}
\begin{proof}
By Proposition \ref{P:comparison1}, (b)$\Longrightarrow$(a). If 
$\operatorname{rr}(X) \leq n$, let  $\epsilon > 0$ and  $(x_{1},\dots ,x_{n+1})$ be an $(n+1)$-tuple of self-adjoint elements in $X$. By Proposition \ref{P:rjs}, there exists an $(n+1)$-tuple $(y_{1},\dots ,y_{n+1})$ of self-adjoint elements of $X$ such that $\| x_{k}-y_{k}\| < \epsilon$ for each $k = 1,\dots ,n+1$ and $\mathbf{0} \not\in \widetilde{\sigma}_{n+1}(y_{1},\dots ,y_{n+1})$. Then Lemma
\ref{L:bremma2} implies that  $(y_{1},\dots ,y_{n+1})$ is $K^2$-unessential. Thus, $\operatorname{br}_{K^2}(X) \leq n$.
\end{proof}

\vspace{0.4in}

%%%%%%%%%%%%%%%%%%%%%%%%%%%%%%%%%%%%%%%%%%%%%%%%%%%%%%%

\subsection{Further properties of bounded rank}\label{SS:fpbr}
As it has been already noted in the Introduction, we can show now that all differences between ``invertible" and $K$-unessential commuting tuples  dissappear when $K\leq 1$. 

\begin{pro}\label{P:commutting}
Let $(y_{1},\dots ,y_{m})$ be a commuting $m$-tuple of self-adjoint elements of the unital $C^{\ast}$-algebra $X$. If $\sum_{i=1}^{m}y_{i}^{2}$ is invertible, then $(y_{1},\dots ,y_{m})$ is $K$-unessential for any positive $K\leq 1$.
\end{pro}
\begin{proof}
It suffuces to prove that $(y_{1},\dots ,y_{m})$ is $1$-unessential.
Take a rational number $\delta >0$  and consider the $C^{\ast}$-subalgebra $Y$ of $X$ generated by the given commuting $m$-tuple $(y_{1},\dots ,y_{m})$ and the unit element. By Gelfand duality there is a compact space $T$
such that $Y$ is isometrically isomorphic to $C(T)$. Consequently we can identify elements of our $m$-tuple with real-valued functions defined on $T$. Consider the diagonal product $y = \triangle \{ y_{i} \colon i =1,\dots ,n+1\} \colon T \to \mathbb{R}^{m}$ defined by $y(t) = \{ y_{1}(t),\dots ,y_{m}(t)\} \in \mathbb{R}^{m}$ for each $t \in T$. Since $\sum_{i=1}^{m}y_{i}^{2}$ is invertible, the image $y(T)$ does not contain the origin $\mathbf{0} = \{ 0,\dots ,0\} \in \mathbb{R}^{m}$. Let $O(\mathbf{0},\delta )$ denote the open ball at $\mathbf{0}$ of radius $\delta$. Consider a retraction

\[ r_{\delta} \colon (\mathbb{R}^{m} \setminus \{\mathbf{0}\} ) \to (\mathbb{R}^{m} \setminus O(\mathbf{0},\delta )) .\]

Since $\mathbf{0} \not\in y(T)$,  the composition $r_{\delta}\circ y \colon T \to \mathbb{R}^{m} \setminus O(\mathbf{0},\delta )$ is well defined. Let $z_{i} = \pi_{i} \circ r_{\delta}\circ y \colon T \to \mathbb{R}$, where $\pi_{i} \colon \mathbb{R}^{m} \to \mathbb{R}_{i} = \mathbb{R}$ denotes the projection onto the $i$-th coordinate. Note that $\| y_{i}-z_{i}\| \leq \delta$ for each $i = 1,\dots , m$. 

Take a point $t \in T$. Since $( z_{1}(t),\dots ,z_{m}(t)) =  r_{\delta}(y_{1}(t), \dots ,y_{m}(t) )\not\in O(\mathbf{0},\delta )$, it follows that $z_{j}(t) \not\in (-\delta ,\delta )$ for some $j = 1,\dots ,m$. 
Consequently $\left( \sum_{i=1}^{m}z_{i}^{2}\right) (t) = \sum_{i=1}^{m}z_{i}^{2}(t) \geq z_{j}^{2}(t) \geq \delta^{2}.$ By \cite[Example 1.2.1]{murphy}, this means that $\sigma\left( \sum_{i=1}^{m}z_{i}^{2}\right) \subseteq [\delta^2 ,\infty)$. This, in turn, implies that
$\left\|\left( \sum_{i=1}^{m}z_{i}^{2} \right)^{-1}\right\| \leq \frac{1}{\delta^2}$. 
\end{proof}

\begin{cor}\label{C:rrcommutative}
Let $X$ be a commutative unital $C^{\ast}$-algebra and $0<K\leq 1$. Then $\operatorname{br}_{K}(X) = \operatorname{rr}(X) = \dim \Omega (X)$, where $\Omega (X)$ is the spectrum of $X$.
\end{cor}
\begin{proof}
By Proposition \ref{P:comparison1}, $\operatorname{rr}(X) \leq \operatorname{br}_{K}(X)$. The opposite inequality $\operatorname{br}_{K}(X) \leq \operatorname{rr}(X)$ follows from Proposition \ref{P:commutting}. The remaining part is well known (see \cite[Proposition 1.1]{brownped91}.
\end{proof}

\begin{cor}\label{C:rrzero}
Let $X$ be a unital $C^{\ast}$-algebra and $0<K\leq 1$. Then $\operatorname{br}_{K}(X) = 0$ if and only if  $\operatorname{rr}(X) =0$.
\end{cor}
\begin{proof}
According to Proposition \ref{P:comparison1}, $ \operatorname{br}_{K}(X)=0$ yields $\operatorname{rr}(X)=0$. Conversely, if $\operatorname{rr}(X) =0$, then, by Proposition \ref{P:commutting}, $\operatorname{br}_{K}(X)=0$.
\end{proof}

Proofs of the following lemma and its corollary are straightforward.
\begin{lem}\label{L:image}
Let $K >0$ and $p\colon X \to Y$ be a surjective $\ast$-homomorphism of unital $C^{\ast}$-algebras. If $(x_{1},\dots ,x_{m})$ is a $K$-unessential $m$-tuple of self-adjoint elements in $X$ then the $m$-tuple $(p(x_{1}),\dots ,p(x_{m}))$ is $K$-unessential in $Y$.
\end{lem}
\begin{proof}
Since $(x_{1},\dots ,x_{m})$ is a $K$-unessential $m$-tuple of self-adjoint elements in $X$,
for any rational $\delta>0$  there exists an $m$-tuple $(y_{1},\dots ,y_{m})$ of self-adjoint elements in $X$ such that
\begin{itemize}
\item[(i)]
$\| x_{k} -y_{k}\| \leq \delta$ for each $k=1,\dots ,m$.
\item[(ii)]
The element $\sum_{k=1}^{m}z_{k}^{2}$ is invertible and $\left\| \left(\sum_{k=1}^{m}z_{k}^{2}\right)^{-1}\right\| \leq \frac{1}{K\cdot \delta}$.
\end{itemize}
Clearly $\| p(x_{k}) -p(y_{k})\| \leq \| x_{k}-y_{k}\| \leq \delta$ for each $k=1,\dots ,m$ and 

\begin{multline*}
\left\| \left(\sum_{k=1}^{m}\left[p(y_{k})\right]^{2}\right)^{-1}\right\| =  \left\|\left( \sum_{k=1}^{m}p(y_{k}^{2})\right)^{-1}\right\| = \left\|\left( p\left(\sum_{k=1}^{m}y_{k}^{2}\right)\right)^{-1}\right\| =\\ \left\| p\left( \left(\sum_{k=1}^{m}y_{k}^{2}\right)^{-1}\right)\right\| \leq \left\| \left(\sum_{k=1}^{m}y_{k}^{2}\right)^{-1}\right\| \leq \frac{1}{K\cdot \delta^2}
\end{multline*}

\noindent which shows that $(p(x_{1}),\dots ,p(x_{m}))$ is indeed $K$-unessential.
\end{proof}

As a corollary we obtain the following useful statement.

\begin{pro}\label{P:image}
Let $K >0$ and $p \colon X \to Y$ be a surjective $\ast$-homomorphism of unital $C^{\ast}$-algebras. Then $\operatorname{br}_{K}(Y) \leq \operatorname{br}_{K}(X)$.
\end{pro}

The following result is actually one of the main reasons of introducing the concept of bounded rank (see the corresponding discussion in the Introduction). It is used in the proof of Theorem \ref{T:main}.
\begin{pro}\label{P:product}
Let $K >0$ and $\{ X_{t} \colon t \in T\}$ be a family of unital $C^{\ast}$-algebras such that
$\operatorname{br}_{K}(X_{t}) \leq n$ for each $t \in T$. Then
$\operatorname{br}_{K}\left(\prod\{ X_{t} \colon t \in T\} \right) \leq n$.
\end{pro}
\begin{proof}
Let $(x_{1},\dots ,x_{n+1})$ be an $(n+1)$-tuple of self-adjoint elements of the product $X = \prod\{ X_{t} \colon t \in T\}$, where $x_{k} = \{ x_{k}^{t} \colon t \in T\}$ for each $k = 1,\dots ,n+1$, and let  $\epsilon > 0$. Our goal is to find a $K$-unessential $(n+1)$-tuple $(y_{1},\dots ,y_{n+1})$ in $X$ such that $\| x_{k} -y_{k}\| < \epsilon$.
For a given $t \in T$ consider the $(n+1)$-tuple $(x_{1}^{t},\dots ,x_{n+1}^{t})$ of self-adjoint elements in $X_{t}$. Since $\operatorname{br}_{K}(X_{t}) \leq n$, there exists a $K$-unessential $(n+1)$-tuple $(y^{t}_{1}, \dots ,y^{t}_{n+1})$ in $X_{t}$ such that $\| x_{k}^{t} - y_{k}^{t}\|_{X_{t}} < \frac{\epsilon}{2}$ for each $k = 1,\dots ,n+1$. Consider the $(n+1)$-tuple $(y_{1}, \dots ,y_{n+1})$, where $y_{k} = \{ y_{k}^{t} \colon t \in T\}$ for each $k = 1,\dots ,n+1$. Note that $y_{k} \in X$ for each $k=1,\dots ,n+1$. Indeed
$\| y_{k}^{t}\|_{X_{t}} \leq \| y_{k}^{t}-x_{k}^{t}\|_{X_{t}} +\| x_{k}\|_{X_{t}} \leq \frac{\epsilon}{2}+\sup\{ \| x_{k}^{t}\|_{X_{t}} \colon t \in T\}$ and $\sup\{ \| y_{k}^{t}\|_{X_{t}} \colon t \in T\} \leq \frac{\epsilon}{2} +\sup\{ \| x_{k}^{t}\|_{X_{t}} \colon t \in T\} < \infty$. Also note that $\| x_{k} -y_{k}\| = \sup \{ \| x_{k}^{t}-y_{k}^{t}\|_{X_{t}} \colon t \in T\} \leq \frac{\epsilon}{2}<\epsilon$.
It only remains to show that the $(n+1)$-tuple $(y_{1}, \dots ,y_{n+1})$ is $K$-unessential in $X$. Indeed, let $\delta > 0$ be rational. Since the $(n+1)$-tuple $(y_{1}^{t},\dots ,y_{n+1}^{t})$ is $K$-unessential in $X_{t}$, there exists an $(n+1)$-tuple $(z_{1}^{t},\dots ,z_{n+1}^{t})$ in $X_{t}$ such that
\begin{itemize}
\item[(a)$_{t}$]
$\| y_{k}^{t}-z_{k}^{t}\|_{X_{t}} \leq \delta$ for each $k=1,\dots ,n+1$.
\item[(b)$_{t}$]
$\left\| \left( \sum_{k=1}^{n+1}\left( z_{k}^{t}\right)^{2}\right)^{-1}\right\| \leq \frac{1}{K\cdot \delta^2}$.
\end{itemize}

Next, consider the $(n+1)$-tuple $(z_{1},\dots ,z_{n+1})$, where $z_{k} = \{ z_{k}^{t} \colon t \in T\}$ for each $k =1,\dots ,n+1$. As above, $z_{k} \in X$ and obviously
\[ \| y_{k}-z_{k}\| = \sup\{ \| y_{k}^{t}-z_{k}^{t}\|_{X_{t}}\colon t \in T \} \leq \delta , k =1,\dots ,n+1 .\]

Finally, note that 
\begin{multline*}
\left\|\left(\sum_{k=1}^{n+1}z_{k}^{2}\right)^{-1}\right\| = \left\| \left\{ \left( \sum_{k=1}^{n+1}\left(z_{k}^{t}\right)^{2}\right)^{-1} \colon t \in T \right\}\right\| =\\ \sup\left\{ \left\| \left( \sum_{k=1}^{n+1}\left( z_{k}^{t}\right)^{2}\right)^{-1} \right\|_{X_{t}} \colon t \in T \right\} \leq \frac{1}{K\cdot \delta^{2}}
\end{multline*}

\noindent as required. 
\end{proof}
%%%%%%%%%%%%%%%%%%%%%%%%%%%%%%%%%%%%

\vspace{0.4in}

%%%%%%%%%%%%%%%%%%%%%%%%%%%%%%%%%%%
%%%%%%%%%%%%%%%%%%%%%%%%%%%%%%%%

\section{Spectral decompositions of unital $C^{\ast}$-algebras of bounded rank $n$}\label{S:decomposition}

In this Section we investigate the behavior of the bounded rank with respect to direct systems and then use this to prove the existence of universal elements in ${\mathcal B}{\mathcal R}_{n}^{K}$.

\begin{pro}\label{P:inverse}
Let $K > 0$ and let $X =  \varinjlim{\mathcal S}$, where ${\mathcal S} = \{ X_{\alpha}, i_{\alpha}^{\beta}, A\}$ is a direct system consisting of unital $C^{\ast}$-algebras and unital $\ast$-inclusions.
If $\operatorname{br}_{K}(X_{\alpha}) \leq n$ for each $\alpha \in A$, then $\operatorname{br}_{K}(X) \leq n$.
\end{pro}
\begin{proof}  Since $\bigcup\{X_{\alpha} \colon \alpha \in A\}$ is dense in $X$ (we identify each $i_{\alpha}(X_{\alpha})$ with $X_{\alpha}$),  
for any $(n+1)$-tuple $(x_{1},\dots ,x_{n+1})$ of self-adjoint elements of $X$ and $\epsilon>0$, we can first find $\alpha\in A$ and an $(n+1)$-tuple $(y_{1},\dots ,y_{n+1})$ of self-adjoint elements of $X_{\alpha}$ with $\| x_{k}^{t}-y_{k}^{t}\|<\frac{\epsilon}{2}$ for each $k=1,\dots ,n+1$. Then, using that $\operatorname{br}_{K}(X_{\alpha}) \leq n$, there exists a $K$-unessential 
$(n+1)$-tuple $(z_{1},\dots ,z_{n+1})$ in $X_{\alpha}$ such that  $\| y_{k}^{t}-z_{k}^{t}\|<\frac{\epsilon}{2}$,  $k=1,\dots ,n+1$. Obviously, $\| x_{k}^{t}-z_{k}^{t}\|<\epsilon$ for each $k=1,\dots ,n+1$ and $(z_{1},\dots ,z_{n+1})$, being $K$-unessential in
$X_{\alpha}$, is $K$-unessential in $X$.
\end{proof}

The following decomposition theorem is known to be true for $n=0$ \cite{chi001}. Proof below works for the real rank as well.

\begin{pro}\label{P:rrzero}
Let $K >0$. The following conditions
are equivalent for any unital $C^{\ast}$-algebra $X$:
\begin{enumerate}
\item
$\operatorname{br}_{K}(X) \leq n$.
\item
$X$ can be represented as the direct limit of a direct
$C_{\omega}^{\ast}$-system
$\{ X_{\alpha}, i_{\alpha}^{\beta}, A\}$ satisfying the following properties:
\begin{itemize}
\item[(a)]
The indexing set $A$ is cofinal and $\omega$-closed in the
$\omega$-complete set $\exp_{\omega}Y$
for some (any) countable dense subset $Y$ of $X$.
\item[(b)]
$X_{\alpha}$ is a $C^{\ast}$-subalgebra of $X$ such that
$\operatorname{br}_{K}(X_{\alpha}) \leq n$, $\alpha \in A$.
\end{itemize}
\end{enumerate}
\end{pro}
\begin{proof}
The implication $(2) \Longrightarrow (1)$ follows from
Proposition \ref{P:inverse}. 

In order to prove the implication $(1) \Longrightarrow (2)$
we first consider a direct $C_{\omega}^{\ast}$-system
${\mathcal S}_{X} = \{ X_{\alpha}, i_{\alpha}^{\beta}, A\}$
with properties indicated in Proposition \ref{P:exists}. Next
consider the following relation $L \subseteq A^{2}$:

\begin{multline*}
 L = \left\{ (\alpha ,\beta ) \in A^{2} \colon \alpha \leq \beta
\; \text{for each $(n+1)$-tuple}\;\; (x_{1},\dots ,x_{n+1})\; \text{in}\;  \left( X_{\alpha}\right)_{sa}\;\text{and}\right.\\
\;\;\;\;\;\;\;\;\;\;\;\text{each}\; \epsilon > 0\;\;\text{there exists a $K$-unessential $(n+1)$-tuple}\; (y_{1},\dots ,y_{n+1}) \;\text{in}\;  X_{\beta}\\
\;\;\;\;\;\;\;\;\;\;\;\left.\text{such that}\; \left\| x_{k} - y_{k} \right\| < \epsilon\; \text{for each}\; k = 1,\dots ,n+1 \right\}.
\;\;\;\;\;\;\;\;\;\;\;\;\;\;\;\;\;\;\;\;\;\;\;\;\;\;\;\;\;\;\;\;\;\;\;\;\;\;\;\;\;\;\;\;\;\;\;\;\;\;\;\;\;\;\;\;\;\;\;\;\;\;\;\;\; 
\end{multline*}

\noindent Let us verify the conditions of \cite[Proposition 1.1.29]{book}. 

{\em Existence}. Let $\alpha \in A$. We need to show that there exists $\beta \in A$ such that $(\alpha ,\beta ) \in L$. For an $(n+1)$-tuple
$\bar{x} = (x_{1},\dots ,x_{n+1}) \in \left( X_{\alpha}\right)_{sa}^{n+1}$ and for a positive integer $m \in {\mathbb N}$ 
first we prove the following assertion:

\begin{itemize}
\item[$(\ast )_{\left(\alpha ,\bar{x},m\right)}$]
There exist an index $\beta\left(\alpha ,\bar{x},m\right)
\in A$, $\beta \left(\alpha,\bar{x},m\right)\geq \alpha$ and a $K$-unessential $(n+1)$-tuple
$\displaystyle \bar{y} = \left( y_{1},\dots ,y_{n+1}\right)$ in $ X_{\beta\left(\alpha ,\bar{x},m\right)}$
such that $\|  x_{k} - y_{k}\| < \frac{1}{m}$ for each $k = 1, \dots ,n+1$.
\end{itemize}

{\em Proof of $(\ast )_{\left(\alpha ,\bar{x},m\right)}$}.
Since $\operatorname{br}_{K}(X) \leq n$, there exists
a $K$-unessential $(n+1)$-tuple $\bar{y} =  \left( y_{1},\dots ,y_{n+1}\right)$ in $ X$  such that
$\|  x_{k} - y_{k}\| < \frac{1}{m}$ for each $k = 1, \dots ,n+1$. Since ${\mathcal S}_{X}$
is a direct $C^{\ast}_{\omega}$-system,
it follows from Lemma \ref{L:strong} that for each $k = 1,\dots ,n+1$, there exists an index
$\beta (\alpha ,\bar{x},m,k) \in A$ such that $y_{k} \in X_{\beta(\alpha ,\bar{x},m,k)}$.
Since $A$ is a $\omega$-complete directed set there exists, by \cite[Corollary 1.1.28]{book}, an index $\gamma \in A$ such that $\alpha ,\beta (\alpha ,\bar{x},m,k) \leq \gamma$ for each $k = 1,\dots ,n+1$.
This obviously implies that
\[ X_{\alpha} \cup \left\{ y_{1}, \dots ,y_{n+1}\right\} \subseteq X_{\alpha} \cup
\bigcup_{k=1}^{n+1} X_{\beta (\alpha ,\bar{x},m,k)} \subseteq X_{\gamma}.\] 

\noindent Note that $\bar{y}$ is $K$-unessential in $X$, but it might not be $K$-unessential in $X_{\gamma}$.  We are going to find $\beta\left(\alpha ,\bar{x},m\right) \in A$ such that $\bar{y}$ is $K$-unessential in $X_{\beta\left(\alpha ,\bar{x},m\right)}$. To this end, for every rational $\delta>0$ fix an $(n+1)$-tuple $z_{\delta}=(z_{1},\dots ,z_{n+1})$ of self-adjoint elements of $X$ such that $\| y_{k}-z_{k}\| \leq\delta$ for each $k = 1,\dots ,n+1$, 
the element $z_{\delta}^2=\sum_{k=1}^{m}z_{k}^{2}$ is invertible in $X$ and $\left\| \left(z_{\delta}^{2}\right)^{-1}\right\| \leq \frac{1}{K\cdot \delta^{2}}$. Next, we choose $\gamma(\delta)\in A$ with $X_{\gamma(\delta)}$ containing the set $\left\{z_{1},\dots ,z_{n+1}\right\}$ and $\gamma (\delta)\geq\gamma$.  Finally, by \cite[Corollary 1.1.28]{book},
there exists an element $\beta (\alpha ,\bar{x},m) \in A$ such that
$\beta (\alpha ,\bar{x},m) \geq \gamma(\delta)$ for every rational $\delta$. It follows from our construction that $\bar{y}$ is a $K$-unessential $(n+1)$-tuple in $X_{\beta(\alpha,\bar{x},m)}$.  This finishes the proof of
$(\ast )_{\left(\alpha , x , m\right)}$. 

For a given $(n+1)$-tuple $\bar{x} = (x_{1},\dots ,x_{n+1}) \in \left( X_{\alpha}\right)_{sa}^{n+1}$
consider indeces
$\beta\left(\alpha ,\bar{x},m\right) \in A$, 
satisfying conditions
$(\ast )_{\left(\alpha ,\bar{x},m\right)}$,
$m = 1,2,\dots$.  Applying \cite[Corollary 1.1.28]{book} once more,
there exists an element $\beta (\alpha ,\bar{x}) \in A$ such that
$\beta (\alpha ,\bar{x}) \geq \beta \left(\alpha,\bar{x} ,m \right)$
for each $m = 1,2,\dots$. Obviously
$X_{\alpha} \subseteq X_{\beta (\alpha ,\bar{x})}$.
Note also that
\begin{itemize}
\item[$(\ast )_{(\alpha ,\bar{x})}$]
For any $\epsilon > 0$ there is a $K$-unessential $(n+1)$-tuple $\bar{y} =\left( y^{\left(\alpha ,\bar{x} \right)}_{1},\dots ,y^{\left(\alpha ,\bar{x}\right)}_{n+1}\right)$ in $ X_{\beta (\alpha ,\bar{x})}$
such that $\left\|  x_{k} - y_{k}^{\left(\alpha ,\bar{x} \right)}\right\| < \epsilon$ for each $k = 1,\dots ,n+1$.
\end{itemize}

Since ${\mathcal S}_{X}$ is a direct $C_{\omega}^{\ast}$-system, it
follows (Definition \ref{D:smooth}(b)) that $X_{\alpha}$ is separable.
So, there exists a countable dense subset
$Y_{\alpha}$ of $\left( X_{\alpha}\right)_{sa}$. For each
$\bar{y} = (y_{1},\dots ,y_{n+1}) \in  Y_{\alpha}^{n+1}$
 consider an index
$\beta (\alpha ,\bar{y}) \in A$ satisfying condition
$(\ast )_{(\alpha ,\bar{y})}$. Since $A$ is a $\omega$-complete set and since $\left| Y_{\alpha}^{n+1}\right| \leq \omega$,
we conclude, by Corollary \cite[Corollary 1.1.28]{book}, that
there exists an index $\beta = \beta (\alpha) \in A$ such that
$\beta  \geq \beta (\alpha ,\bar{y})$ for each
$\bar{y} \in Y_{\alpha}^{n+1}$.

We claim that $(\alpha ,\beta ) \in L$. Indeed, let $\epsilon > 0$ and
$\bar{x} = (x_{1},\dots ,x_{n+1}) \in \left( X_{\alpha}\right)_{sa}^{n+1}$.
Since, by the above construction,
$Y_{\alpha}$ is dense in
$\left( X_{\alpha}\right)_{sa}$, there exists an $(n+1)$-tuple
$\tilde{y} = (\tilde{y}_{1}\dots ,\tilde{y}_{n+1}) \in Y_{\alpha}^{n+1}$ such that 
$\| \tilde{y}_{k} - x_{k}\| \leq \frac{\epsilon}{2}$ for each $k = 1,\dots , n+1$. By condition
$(\ast )_{(\alpha ,\tilde{y})}$, there exists a $K$-unessential $(n+1)$-tuple 
$\bar{y} = (y_{1}^{(\alpha ,\tilde{y})},\dots ,y_{n+1}^{(\alpha ,\tilde{y})}) $ in $ X_{\beta(\alpha ,\tilde{y})}$
such that $\| \tilde{y}_{k} -y_{k}^{(\alpha ,\tilde{y})}\| < \frac{\epsilon}{2}$
for each $k = 1,\dots ,n+1$.
Since
$\beta =\beta(\alpha) \geq \beta (\alpha ,\tilde{y})$, it follows that
$X_{\beta (\alpha ,\tilde{y})} \subseteq X_{\beta (\alpha )}$.
This guarantees that $\bar{y} \in \left( X_{\beta (\alpha )}\right)_{sa}^{n+1}$.
It only remains to note that $\| x_{k} - y_{k}^{(\alpha ,\tilde{y})} \| \leq \| x_{k} -\tilde{y}_{k}\| +\| \tilde{y}_{k} - y_{k}^{(\alpha ,\tilde{y})}\| < \epsilon$.
Therefore $(\alpha, \beta ) \in L$.

{\em Majorantness}. Let $(\alpha ,\beta ) \in L$ and 
$\gamma \geq \beta$. We need to show that $(\alpha ,\gamma ) \in L$. Let $\epsilon > 0$ and
$x = (x_{1},\dots ,x_{n+1}) \in \left( X_{\alpha}\right)_{sa}^{n+1}$.
Since $(\alpha ,\beta ) \in L$, there exists a $K$-unessential $(n+1)$-tuple
 $y = (y_{1},\dots ,y_{n+1})$ in  $ X_{\beta}$ such that $\| x_{k}-y_{k}\| < \epsilon$ for each $k = 1,\dots ,n+1$. Since $\gamma \geq \beta$, it follows
that
$\left( X_{\beta}\right)_{sa} \subseteq \left(
X_{\gamma}\right)_{sa}$
which shows that $y \in \left( X_{\gamma}\right)_{sa}^{n+1}$. Moreover, $y$ is  $K$-unessential in $X_{\gamma}$, so
$(\alpha ,\gamma ) \in L$.

{\em $\omega$-closeness}. Suppose that $\{ \alpha_{i} : i \in \omega \}$
is a countable chain of indeces in $A$. Assume also that
$(\alpha_{i}, \beta ) \in L$ for each $i \in\omega$ and some
$\beta \in A$. Our goal is to show that
$(\alpha ,\beta ) \in L$, where $\alpha =
\sup \{\alpha_{i} \colon i \in \omega \}$. Let $\epsilon > 0$ and
$(x_{1},\dots ,x_{n+1}) \in \left( X_{\alpha}\right)_{sa}^{n+1}$.
Since ${\mathcal S}_{X}$ is a direct $C^{\ast}_{\omega}$-system
it follows (Definition \ref{D:smooth}(c)) that
$X_{\alpha}$ is the direct limit of the direct system
generated by $C^{\ast}$-subalgebras
$X_{\alpha_{i}}$, $i \in \omega$, and corresponding inclusion homomorphisms.
Consequently, there exist an index $j \in\omega$ and an $(n+1)$-tuple 
$(x_{1}^{j},\dots ,x_{n+1}^{j}) \in \left( X_{\alpha_{j}}\right)_{sa}^{n+1}$ such that
$\| x_{k} - x_{k}^{j}\| < \frac{\epsilon}{2}$ for each $k = 1,\dots ,n+1$. Since
$(\alpha_{j},\beta ) \in L$, there exists a $K$-unessential $(n+1)$-tuple 
$(y_{1}^{j},\dots ,y_{n+1}^{j})$ in $ X_{\beta}$ such that $\| x_{k}^{j} - y_{k}^{j}\| <\frac{\epsilon}{2}$ for each $k = 1,\dots ,n+1$. Clearly
$\| x_{k}-y_{k}^{j}\| < \epsilon$ for each $k = 1,\dots ,n+1$. This shows that $(\alpha ,\beta ) \in L$.

We are now in a position to apply \cite[Proposition 1.1.29]{book} which
guarantees that
the set $A^{\prime} = \{ \alpha \in A \colon (\alpha ,\alpha ) \in L\}$ is cofinal
and $\omega$-closed in $A$.
Note that 
$(\alpha ,\alpha ) \in L$ precisely when for each $\epsilon > 0$ and for
each $(n+1)$-tuple $(x_{1},\dots ,x_{n+1}) \in \left( X_{\alpha}\right)_{sa}^{n+1}$ there exists a $K$-unessential $(n+1)$-tuple $(y_{1},\dots ,y_{n+1}) $ in  $ X_{\alpha}$ such
that $\| x_{k}-y_{k}\| < \epsilon$ for each $k = 1,\dots ,n+1$. 
This obviously means that the direct $C_{\omega}^{\ast}$-system
$\varinjlim{\mathcal S}_{X}^{\prime} =
\{ X_{\alpha}, i_{\alpha}^{\beta}, A^{\prime}\}$ consists of
$C^{\ast}$-subalgebras of $X$ of bounded rank (with respect to $K$) at most $n$. Clearly
$\varinjlim{\mathcal S}_{X}^{\prime} = X$.
Proof is completed.
\end{proof}

The following statement immediately follows form Proposition \ref{P:rrzero}. 
\begin{cor}\label{C:separable}
Let $K >0$. Every countable subset of a unital $C^{\ast}$-algebra $X$ with $\operatorname{br}_{K}(X) \leq n$  
is contained in a unital separable $C^{\ast}$-subagebra $X_{0}$ such that $\operatorname{br}_{K}(X_{0}) \leq n$.
\end{cor}

Next, for any $K > 0$, we construct a universal separable unital
$C^{\ast}$-algebra $Z_{n}^{K}$ of bounded rank $n$. Universal in the
sense that any other
separable unital $C^{\ast}$-algebra with bounded rank $\leq n$ is its
quotient.

\begin{thm}\label{T:main}
Let $K >0$. The class ${\mathcal B}{\mathcal R}_{n}^{K}$ of all separable unital
$C^{\ast}$-algebras with $\operatorname{br}_{K} \leq n$ contains an universal element $Z_{n}^{K}$.
More formally,
there exists a ${\mathcal B}{\mathcal R}_{n}^{K}$-invertible
unital $\ast$-homomorphism
$p \colon C^{\ast}\left( {\mathbb F}_{\infty}\right)
\to Z_{n}^{K}$, where $Z_{n}^{K}$ is a separable unital
$C^{\ast}$-algebra such that $\operatorname{br}_{K}\left( Z_{n}^{K}\right) = n$.
\end{thm}
\begin{proof}
Let
${\mathcal A} =
\{ f_{t} \colon C^{\ast}\left( {\mathbb F}_{\infty}\right)
\to X_{t}, t \in T\}$ denote the set of all unital
$\ast$-ho\-mo\-mor\-phisms, defined on
$C^{\ast}\left( {\mathbb F}_{\infty}\right)$,
such that $\operatorname{br}_{K}(X_{t}) \leq n$. Next, consider the product
$\prod\{ X_{t} \colon t \in T\}$. Since
$\operatorname{br}_{K}(X_{t}) \leq n$ for each $t \in T$,
it follows from Proposition \ref{P:product} that
$\operatorname{br}_{K}\left( \prod\{ X_{t} \colon t \in T\}\right) \leq n$.
The $\ast$-homomorphisms $f_{t}$, $t \in T$, define the
unital $\ast$-homomorphism
$f \colon C^{\ast}\left( {\mathbb F}_{\infty}\right)
\to \prod\{ X_{t} \colon t \in T\}$ such that $\pi_{t} \circ f = f_{t}$ for each $t \in T$ (here
$\pi_{t} \colon \prod\{ X_{t} \colon t \in T\} \to X_{t}$
denotes
the corresponding canonical projection $\ast$-homomorphism).
By Proposition \ref{P:rrzero},
$\prod\{ X_{t} \colon t \in T\}$ can be
represented as the limit of the $C^{\ast}_{\omega}$-system
${\mathcal S} = \{ C_{\alpha}, i_{\alpha}^{\beta}, A\}$ such
that $C_{\alpha}$ is a separable unital $C^{\ast}$-algebra
with $\operatorname{br}_{K}(C_{\alpha})\leq n$ for each $\alpha \in A$. Suppressing injective unital
$\ast$-homomorphisms $i_{\alpha}^{\beta} \colon C_{\alpha} \to C_{\beta}$, we can,
for notational simplicity, assume that $C_{\alpha}$'s are unital
$C^{\ast}$-subalgebras of $\prod\{ X_{t} \colon t \in T\}$.
Let $\{ a_{k} \colon k \in \omega\}$ be a countable dense subset of
$C^{\ast}\left( {\mathbb F}_{\infty}\right)$. By Lemma \ref{L:strong},
for each $k \in \omega$ there exists an index $\alpha_{k} \in A$ such that
$f(a_{k}) \in C_{\alpha_{k}}$. By \cite[Corollary 1.1.28]{book}, there exists an index
$\alpha_0 \in A$ such that $\alpha_0 \geq \alpha_{k}$ for each $k \in \omega$.
Then 
$f(a_{k}) \in C_{\alpha_{k}} \subseteq C_{\alpha_{0}}$ for each $k \in \omega$ (see also Corollary \ref{C:separable}).
This observation coupled with the continuity of $f$ guarantees that
$f\left( C^{\ast}\left( {\mathbb F}_{\infty}\right)\right) =
f\left(\operatorname{cl}\left \{ a_{k} \colon k \in \omega\right\}\right)
\subseteq \operatorname{cl}\left\{ f\left(\{ a_{k} \colon k \in \omega\}
\right)\right\} \subseteq \operatorname{cl}C_{\alpha_{0}} = C_{\alpha_{0}}$.

Let $Z_{n}^{K} = C_{\alpha_{0}}$ and $p$ denote the unital $\ast$-homomorphism $f$
considered as the homomorphism of
$C^{\ast}\left( {\mathbb F}_{\infty}\right)$ into $Z_{n}^{K}$.
Note that $f = i \circ p$, where
$i \colon Z_{n}^{K} = C_{\alpha_{0}}
\hookrightarrow \prod\{ X_{t} \colon t \in T\}$ stands
for the inclusion.

By construction, $\operatorname{br}_{K}(Z_{n}^{K}) \leq n$. Let us show that
$\displaystyle p \colon C^{\ast}\left( {\mathbb F}_{\infty}\right) \to Z_{n}^{K}$
is ${\mathcal B}{\mathcal R}_{n}^{K}$-invertible in the sense of Introduction.
In our situation, for any unital $\ast$-homo\-mor\-phism
$g \colon C^{\ast}\left( {\mathbb F}_{\infty}\right) \to X$, where $X$ is a
separable unital $C^{\ast}$-algebra with $\operatorname{br}_{K}(X)\leq n$, we
need to establish the existence of a unital $\ast$-homomorphism
$h \colon Z_{n}^{K} \to X$ such that $g = h\circ p$. Indeed, by definition
of the set ${\mathcal A}$, we conclude that $g = f_{t}$ for some index $t \in T$
(in particular, $X = X_{t}$ for the same index $t \in T$).
Next observe that $g = f_{t} = \pi_{t} \circ f = \pi_{t} \circ i \circ p$.
This allows us to define the required unital $\ast$-homo\-mor\-phism
$h \colon Z_{n}^{K} \to X$ as the composition $h = \pi_{t} \circ i$. Hence,
$p$ is ${\mathcal B}{\mathcal R}_{n}^{K}$-invertible which yields
 (see Introduction) the universality of $Z_{n}^{K}$.
\end{proof}

%%%%%%%%%%%%%%%%%%%%%%%%%%%%%%%%%%%%%%%%%%%%%%%%%%%%%%%%%%%%%%%%

%%%%%%%%%%%%%%%%%%%%%%%%%%%%%%%%%%
%%%%%%%%%%%%%%%%%%%%%%%%%%%%%%%%%%%%%%%

\vspace{0.4in}

\section{Concluding remarks}\label{S:concluding}

It would be interesting to find conditions under which bounded and real ranks coincide. Some of them, formulated in terms of joint spectra, have been discussed in Section \ref{S:rjs}. Let us consider one more, of a somewhat different nature. Let us say that a unital $C^{\ast}$-algebra $X$ has  an open $m$-squaring map if the map $\alpha_{m} \colon \left( X_{sa}\right)^{m} \to X_{+}$, defined by letting $\alpha_{m}(x_{1},\dots ,x_{m}) = \sum_{k=1}^{m}x_{k}^{2}$, is open.

\begin{pro}\label{P:square}
If the unital $C^{\ast}$-algebras $X$ has an open $(n+1)$-squaring map, then
$\operatorname{rr}(X) \leq n$.
\end{pro}
\begin{proof}
Let $(x_{1},\dots ,x_{n+1})$ be an $(n+1)$-tuple of self-adjoint elements of $X$ and $\epsilon > 0$. Clearly $\sigma (x) \subseteq [0,\infty )$, where $x = \alpha_{n+1}(x_{1},\dots ,x_{n+1}) = \sum_{k=1}^{n+1}x_{k}^{2}$. Since $\alpha_{n+1}$ is open, there exists $\delta > 0$ (generally speaking $\delta$ depends on the $(n+1)$-tuple
$(x_{1},\dots , x_{n+1})$) such that for every $y \in X_{+}$ with $\| x-y\| < \delta$ there exists an $(n+1)$-tuple $(y_{1},\dots ,y_{n+1})$ of self-adjoint elements of $X$ with $\alpha_{n+1}(y_{1}, \dots ,y_{n+1}) = y$ and $\| y_{k} -x_{k}\| < \epsilon$ for each $k = 1,\dots ,n+1$. Next, consider a function $f \colon \sigma (x) \to ( 0,\infty )$ such that $\| f(t) -t\| < \delta$, $t \in \sigma (x)$. Functional calculus supplies the element $z \in C^{\ast}(x)$ corresponding to $f$ such that $\|z -x\| < \delta$. Since $f \in C(\sigma (x))$ is positive (because $f(\sigma (x)) \subseteq (0,\infty )$),  $z$ is also positive. By the Spectral Mapping Theorem, $\sigma (z) = f(\sigma (x)) \subseteq (0,\infty )$. Consequently, $z$ is invertible. Choice of $\delta$ guarantees the existence of an $(n+1)$-tuple  $(z_{1},\dots ,z_{n+1})$ of self-adjoint elements such that $\alpha_{n+1}(z_{1},\dots z_{n+1}) = z$ and $\| z_{k} -x_{k}\| < \epsilon$ for each $k=1,\dots ,n+1$. By Lemma \ref{L:definition}(ii), this means that $\operatorname{rr}(X) \leq n$.
\end{proof}

It would be interesting to understand how stronger the requirement of openness of the $(n+1)$-squaring map is compared to the condition $\operatorname{rr}(X) \leq n$ and when do these two conditions coincide\footnote{This question is discussed in \cite{ckr}.}. 

We are grateful to Vladimir Kisil for information regarding general theory of joint spectra for non-commuting tuples.
%%%%%%%%%%%%%%%%%%%%%%%%%%%%%%%%%%%%%%%%%%%%%%%%%%%%%%%%%%%%%%%

%%%%%%%%%%%%%%%%%%%%%%%%

%%%%%%%%%%%%%%%%%%%%%%%%%%%%%%%%%

%%%%%%%%%%%%%%%%%%%%%%%%%%%%%%%%%%%%%


\begin{thebibliography}{99}




\bibitem{brownped91}
L.~G.~Brown and G.~K.~Pedersen, {\em $C^{\ast}$-algebra
of real rank zero},
J. Functional Anal. {\bf 99} (1991), 131-149.

\bibitem{chi991}
A.~Chigogidze, {\em Uncountable direct systems and a
characterization of
non-separable projective $C^{\ast}$-algebras}, Mat. Stud.
{\bf 12}, \# 2 (1999), 171--204.
 
\bibitem{chi001}
A.~Chigogidze, {\em Universal $C^{\ast}$-algebra of real rank zero}, Infinite Dimensional Analysis, Quantum Probability and Related Topics, {\bf 3} (2000), 445--452.

\bibitem{book}
A.~Chigogidze, {\em Inverse Spectra}, North Holland,
Amsterdam, 1996.

\bibitem{ckr}
A.~Chigogidze, A.~Karasev, M.~R{\o}rdam, {\em Real rank and squaring mappings for unital $C^{\ast}$-algebras}, Preprint math.FA/0201214 (2002).

\bibitem{eng78}
R.~Engelking, {\em Dimension Theory}, PWN, Warsaw, 1978.

\bibitem{hw}
W.~Hurewitz, H.~Wallman, {\em Dimension Theory}, Princeton Univ. Press, 1941.

\bibitem{lin}
H.~Lin, {\em The tracial topological rank of $C^{\ast}$-algebras},
Proc. London Math. Soc. {\bf 83} (2001), 199--234.



\bibitem{murphy}
G.~J.~Murphy, {\em $C^{\ast}$-algebras and Operator Theory},
Academic Press, London, 1990.


\bibitem{murphy2}
G.~J.~Murphy, {\em The analytic rank of a $C^{\ast}$-algebra}, Proc. Amer. Math. Soc. {\bf 115} (1992), 741--746.


\bibitem{nagata}
J.~Nagata, {\em Modern Dimension Theory}, North-Holland, Amsterdam, 1965.


\bibitem{pears}
A.~R.~Pears, {\em Dimension Theory of General Spaces}, Cambridge Univ. Press, London, 1975.

\bibitem{phillips}
N.~C.~Phillips, {\em Simple $C^{\ast}$-algebras with the property (FU)}, Math. Scand. {\bf 69} (1991), 121--151.

\bibitem{rieffel}
M.~A.~Rieffel, {\em Dimension and stable rank in the $K$-theory of $C^{\ast}$-algebras},
Proc. London Math. Soc. {\bf 46} (1983), 301--333.

\bibitem{taylor1}
J.~P.~Taylor, {\em A general framework for a multi-operator functional calculus}, Advances Math. {\bf 9}, (1972), 183--252.

\bibitem{taylor2}
J.~P.~Taylor, {\em Functions of several noncommuting variables}, Bull. Amer. Math. Soc. {\bf 79} (1973), 1--34.

\bibitem{winter}
W.~Winter, {\em Covering dimension for nuclear $C^{\ast}$-algebras}, Preprint math.OA/0107218 (2001).


\bibitem{zelazko}
W.~Zelazko, {\em An axiomatic approach to joint spectra I}, Studia Math. {\bf 64} (1979), 249--261.

\end{thebibliography}
\end{document}